\documentclass[12pt,oneside]{amsart}

      \usepackage{amssymb}
      \usepackage{rotating}

\usepackage{amsfonts}
\usepackage{dsfont}
\usepackage{amsxtra}
\usepackage{graphicx}
\DeclareGraphicsExtensions{.pdf,.png,.jpg}

\usepackage{index}
      \theoremstyle{plain}
      \newtheorem{Theo}{Theorem}[section]
      \newtheorem{Lem}[Theo]{Lemma}
      
      \newtheorem{Prop}[Theo]{Proposition}
      \newtheorem{Cor}[Theo]{Corollary}
      
      \theoremstyle{definition}
      \newtheorem{Defi}[Theo]{Definition}

      \theoremstyle{remark}
      \newtheorem{Rem}[Theo]{Remark}

      \def\ZZ{{\mathds Z}}

      \makeatletter
      \def\@setcopyright{}
      \def\serieslogo@{}
      \makeatother
   
\catcode`\á=\active \def á{\'a}
 \catcode`\Á=\active \def Á{\'A}
  \catcode`\ó=\active \def ó{\'o}
  \catcode`\é=\active \def é{\'e}
  \catcode`\ú=\active \def ú{\'u}
  \catcode`\í=\active \def í{\'{\i}}
  \catcode`\ñ=\active \def ñ{\~n}
  \catcode`\Ñ=\active \def Ñ{\~N}
  \catcode`\¿=\active \def ¿{?`}
  \catcode`\º=\active \def º{$^{\underline{o}}$}
  \catcode`\ª=\active \def ª{$^{\underline{a}}$}
  \catcode`\¡=\active \def ¡{!`}
  \catcode`\â=\active \def â{\^{a}}
  \catcode`\ê=\active \def ê{\^{e}}
  \catcode`\î=\active \def î{\^{\i}}
  \catcode`\ô=\active \def ô{\^{o}}
  \catcode`\û=\active \def û{\^{u}}
  \catcode`\ç=\active \def ç{\c{c}}
  \catcode`\ü=\active \def ü{\"{u}}
  \catcode`\ö=\active \def ö{\"{o}}
 \newenvironment{dem}{\noindent{\it Proof}.}{\qed}

\newcommand{\pro}{\mathbb{P}}
\newcommand{\oo}{\mathcal{O}}
\newcommand{\ooo}{\overline{\mathcal{O}}}
\newcommand{\mm}{\mathfrak{m}}

\newcommand{\Nu}{{ \mbox{{\LARGE $\nu$}} } }

\begin{document}

% First we specify the top matter (author, title, etc).
%
% Note: All of the top matter items are optional and can be omitted.
% But you will probably want to specify at least the author and title
% and perhaps an abstract.

   % author information

   % first author 
   
%  \author{F\'elix Delgado de la Mata}
%  \address{Departamento de Algebra, Geometria y Topologia, Universidad de Valladolid. Prado de la Magdalena, s/n E-47005 Valladolid, Spain}
%  \email{fdelgado@agt.uva.es}
   
%   \author{Evgeny Gorsky}

   % the address where the research was carried out
%   \address{USA}

   % current address, usually not needed because it is the same as the
   % regular address
   %\curraddr{Department of Mathematics, Pennsylvania State University,
   %  University Park, State College PA 16802}

 %  \email{usa@usa.usa}

\author{Julio Jos\'e Moyano-Fern\'andez}
\address{Institut f\"ur Mathematik, Universit\"at Osnabr\"uck. Albrechtstra\ss e 28a, D-49076 Osnabr\"uck, Germany}
\email{jmoyanof@uni-osnabrueck.de}

   % title

   \title[Poincar\'e series and behaviour by projections]{Poincar\'e series for plane curve \\ singularities and their behaviour \\ under projections}

   % Note that the short title for running heads goes in square
   % brackets.  This is optional.  The long title goes in curly
   % braces.  In the long title, line breaks are indicated by \\.

   % abstract (optional)
   \begin{abstract}
     Our purpose is to investigate all defined Poincar\'e series associated with multi-index filtrations and value semigroups of curve singularities---not necessarily complex---with regard to the property of forgetting variables, i.e., by making variables of the series to be $1$. Generalised Poincar\'e series of motivic nature will be also considered.
   \end{abstract}

   % AMS subject classifications (used in AMS journals)

   % AMS keywords (used in AMS journals)

   % acknowledge support, etc
   \thanks{The author was supported by the Deutsche Forschungsgemeinschaft (DFG), and partially by the Spanish Government grant ``Ministerio de Educaci\'on y Ciencia (MEC) 
MTM2007-64704", in cooperation with the European Union in the framework of the founds ``FEDER''}
  % \thanks{We would like to thank our colleagues for their helpful
   %  criticism.?}

   % dedication
   %\dedicatory{Dedicated to Professor Karlheinz Kiyek on the occasion
    %of his $75$th birthday}

   % today's date, or fill in whatever date you prefer
  % \date{\today}

% This ends the top matter information.
% We can now tell LaTeX to display the top matter.

   \maketitle

% Having displayed the top matter, we now proceed to the body of the
% article.

% The body of the article is divided into sections.
% Each section begins with a \section command.
\medskip

\section{Introduction}

The idea of associating a generating function to a ring as a form of characterising properties of the ring goes back to R. Dedekind, as he was able to assimilate Riemann's philosophy in order to define the so-called Dedekind zeta functions for the ring of integers of a number field. The key point was to find a suitable way to measure filtrations of ideals (namely, the norm of an ideal in the ring of integers, since they are all finite). D. Hilbert and E. Noether extended later this argument to commutative algebra by introducing the Hilbert-Poincar\'e series of a module. As soon as the algebraic Geometry ripped in the shade of the fast development of the abstract algebra and the old number-theoretical ideas, similar series were associated to affine rings of algebraic curves, giving rise to zeta functions of curves. A last step---much more recent and less explored---was to deal with singular algebraic curves over finite fields (cf. \cite{Gal}, \cite{Gre}, \cite{St1}). All mentioned zeta functions and Hilbert-Poincar\'e series can be understood as formal power series in one indeterminate. For the case of reduced curve singularities first, and in particular for complex plane curve singularities later, Campillo, Delgado, Kiyek and Gusein-Zade realised that a definition of multivariable Poincar\'e series (i.e. with several indeterminates, one for each irreducible component of the curve) not only makes sense but also yields a finer invariant of the singularity (see \cite{CDK} and \cite{CDG-P}--\cite{cadegu11}). 
\medskip

The topological meaning of changing from one to more variables (for instance by making some variable to take the value $1$) and back was cleared when Campillo, Delgado and Gusein-Zade showed the surprising connection between their Poincar\'e series and the Alexander polynomial for plane complex curve singularities (see \cite{CDG-Alex}), because for this topological invariant, an old paper of Torres had already solved the problem \cite[Theorem 3]{Torres}. (Torres describes even its functional equations, see op.cit. Theorem 2).
\medskip

Nevertheless Poincar\'e series are combinatorial objects definable not only in complex contexts. Our aim in this work is to understand the combinatorics behind the mechanism of passing from one to several variables. We understand that it may be useful, especially after a theorem of Delgado and the author saying that the information in the zeta functions of singular curves over finite fields---objects of number-theoretical nature---is already contained in a generalised Poincar\'e series in the framework of motivic integration (by specialising in the cardinality of the finite field), establishing an interesting bridge between two a priori different objects (see \cite{DM}).  Some related results in the flavour of that theorem were proven later on (e.g.  \cite{MZ}, \cite{Gor}, \cite{St2}, \cite{DM2}).
\medskip

We will keep the following notation throughout the paper. Set $I:=\{1, \ldots , r\}$. For $i \in I$, we will write $\underline{1}_{i}$ for the $i$-th vector of the canonical basis of $\mathds{Z}^r$, and $\underline{1}_J:=\sum_{j \in J} \underline{1}_{i}$ for every $J \subseteq I$. Notice that $\underline{1}_I=\underline{1}:=(1, 1, \ldots ,1)$ and $\underline{1}_{\varnothing}=\underline{0}:=(0, 0, \ldots , 0)$. 
\medskip

\section{Value and extended semigroup}

Let $\oo$ be a one-dimensional Cohen-Macaulay Noetherian local
ring containing a perfect field $k$, with maximal ideal $\mm$. Let
$\mathcal{K}$ be the total ring of fractions of $\oo$ and let
$\overline{\oo}$ be the integral closure of $\oo$ in
$\mathcal{K}$. Assume the degree $\rho:=[\oo/ \mm : k]$ to be
finite. By hypothesis, $\ooo$ is an $\oo$-module of finite length.
\medskip

The integral closure $\ooo$ decomposes as finite intersection of
discrete Manis valuation rings (see \cite{kiyek}): $\ooo = V_1 \cap \ldots \cap V_r$, with
associated discrete Manis valuations $v_1, \ldots, v_r$. For every
$i \in I$,  we will write $\mm_i := \mm (V_i) \cap
\ooo$, where $\mm (V_i)$ is the maximal ideal of $V_i$, and the residue fields will be denoted by $k_i:=V_i/\mm_i$.
\medskip

We will assume that
\[
\oo / \mm = \overline{\oo}/ \mm_i = V_i / \mm (V_i) ~ ~ ~ ~ ~ \textrm{   for every } i \in I.
\]
This hypothesis  implies in particular that $\rho =[k_i:\oo/\mm] =1$ for every $i \in I$, and we will refer to this property by saying that the ring $\oo$ is residually rational (cf. \cite[(4.1)]{CDK}).
\medskip

Thus, given such a ring $\oo$ one can associate the set
\[
S=S_{\oo}:=\{\underline{v}(z) \mid z \in \mathcal{O} \mathrm{~with~} v_i(z)< \infty \mathrm{~for~all~} i \in I \}
\]
which is an additive sub-semigroup of $\mathds{Z}_{\ge 0}^r$: it is said to be the value semigroup associated to the ring $\oo$. Notice that the semigroup $S$ possesses a \emph{conductor}, say $\delta \in S$, i.e., the smallest element $\underline{v} \in S$ such that $\underline{v} + \ZZ^{r}_{\ge 0} \subseteq S$.
\medskip

Let $\underline{v} \in \mathds{Z}$, and $J \subseteq I$. Let us define the following subsets of the value semigroup $S$:

\begin{eqnarray}
\overline{\Delta}^J(\underline{v}) &:=& \{\underline{w} \in \mathbb{N}^r \mid w_i > v_i \mathrm{~ if ~} i \in J \mathrm{~and~}  w_i = v_i \mathrm{~ if ~} i \notin J  \}. \nonumber \\
\Delta^{J} (\underline{v}) &:=& \overline{\Delta}^J(\underline{v}) \cap S. \nonumber \\
\Delta_{J}(\underline{v}) &:= &\Delta^{I \setminus J} (\underline{v}). \nonumber \\
\Delta(\underline{v}) &:=& \bigcup_{i=1}^{r} \Delta^{i}(\underline{n}).  \nonumber 
\end{eqnarray}

A special role will be played by the following elements of $S$ (cf. \cite{D1}):

\begin{Defi} \label{Defi:maximals}
An element $\underline{v} \in S$ is said to be a maximal (of $S$) if $\Delta (\underline{v})= \varnothing$. If, moreover, one has that $\Delta_J (\underline{v})=\varnothing$ for every $J \subsetneq I$, $J \ne \varnothing$, then the element will be called an absolute maximal of $S$. If $\underline{v}$ is maximal and if $\Delta_J (\underline{v}) \ne \varnothing$ for every $J \subseteq I$ such that $\sharp (J) \ge 2$, then $\underline{v}$ will be called a relative maximal of $S$.
\end{Defi}

Let $\oo$ be assumed to have a coefficient field $K$.  Let $\{V_1, \ldots , V_r\}$
be the set of pairwise different Manis valuation rings of
$\mathcal{K}$ belonging to $\oo$.
\medskip

If $K_i$ is a coefficient field of $V_i$ and $t_i$ is an
indeterminate over $K_i$, then one can identify $V_i \cong K_i
[\![t_i]\!]$ and $v_i$ with the order function respect to $t_i$ in
$K_i [\![t_i]\!]$ for every $i \in I$. Hence
\[
\oo \subset K_1 [\![t_1]\!] \cap \ldots \cap K_r[\![t_r]\!] = \ooo.
\]

Since $V_i$ is a module of finite type over the ring $\oo$, the field
extensions $\oo / \mm \hookrightarrow \ooo / \mm_i$ are finite, for
every $i \in I$. Furthermore, as $\oo / \mm$ is
assumed to be perfect, every such a extension is separable and
therefore, for every coefficient field $K$ of $\oo$ there exists a
unique coefficient field $K_i$ of $V_i$ with $K \subset K_i$ which
is isomorphic to $\ooo / \mm_i$ for every $i \in I$.
\medskip

Define the ideals
\[
J(\underline{v}):= \{ g \in \oo \mid \underline{v}(g) \ge
\underline{v}\},
\]
where $\underline{v}(g)=(v_1(g), \ldots , v_r (g))$ and
$\underline{v}(g) \ge \underline{v}$ means $v_i(g) \ge v_i$
for every $i \in I$. For every $\underline{v} \in \mathds{Z}^r$
and every $i \in I$, let us consider the
$\oo$-module defined as
\[
C(\underline{v},i):= J(\underline{v})/J(\underline{v}+\underline{1}_{i}).
\]
Notice that the
$\oo$-module $C(\underline{v},i)$ is annihilated by $\mathfrak{m}$
so that $C(\underline{v},i)$ naturally gets a structure of
$k$-vector space.
\medskip

Let us consider the map
\begin{displaymath}
\begin{array}{lccc}
j_{\underline{v}}: & J(\underline{v}) & \longrightarrow & C(\underline{v},1) \times \ldots \times C(\underline{v},r)   \\
& g & \mapsto & \left ( j_1 (g), \ldots , j_r (g) \right
)=:j_{\underline{v}}(g) .
\end{array}
\end{displaymath}
We can identify $\mathrm{Im}(j_{\underline{v}}) \cong
J(\underline{v})/J(\underline{v}+\underline{1})=:C(\underline{v})$ and define the
set
\[
F_{\underline{v}}:= C(\underline{v}) \cap \left (
(C(\underline{v},1) \setminus \{0\}) \times \ldots \times
(C(\underline{v},r) \setminus \{0\}) \right ).
\]
\begin{Lem}
\[
F_{\underline{v}} = C(\underline{v}) \cap \left (K_1^{\ast}
\times \ldots \times K_r^{\ast} \right ).
\]
\end{Lem}

\proof~It is enough to see that the map $\varphi: C(\underline{v},1) \setminus
\{0 \} \to K_1^{\ast}$ is an isomorphism. Let be $g \in \oo
\setminus \{0 \}$ with $v_1(g)=v_1$. One has that $g = a(t_1)
t_1^{v_1(g)} + \mathrm{~higher~order~terms}$, with $a(t_1) \in
K_1^{\ast}$, thus $\varphi$ can be defined by $g \mapsto a(t_1)$.
\qed

\medskip

In fact, if $\underline{v} \in \mathbb{Z}^r_{\ge 0}$ then one has
\[
F_{\underline{v}}= \left ( J(\underline{v}) /
J(\underline{v}+\underline{1}) \right ) \setminus
\bigcup_{i=1}^{r} \left ( J(\underline{v} + \underline{1}_{i}) / J(\underline{v}+\underline{1}) \right ),
\]
i.e., $F_{\underline{v}}$ is the complement to an arrangement of
vector subspaces in a vector space (but not a vector subspace
itself).

\begin{Defi}
The \emph{extended semigroup}
$\widehat{S}$ is the union of the subspaces $F_{\underline{v}}$
for all $\underline{v} \in \mathds{Z}^r_{\ge 0}$. The spaces
$F_{\underline{v}}$ are called \emph{fibres} of the extended semigroup $\widehat{S}$.
\end{Defi}

The group $K^{\ast}:=K \setminus {0}$ acts freely on
$\mathds{Z}^r \times (K_1^{\ast} \times \ldots \times K_r^{\ast})$
(by multiplication of all coordinates in $K_1^{\ast} \times \ldots
\times K_r^{\ast}$). The corresponding factor space
$\mathds{Z}^{r} \times (K_1^{\ast} \times \ldots \times
K_r^{\ast}) / K^{\ast} = \mathds{Z}^r \times \pro (K_1^{\ast}
\times \ldots \times K_r^{\ast})= \sum_{\underline{v} \in
\mathds{Z}^r} \pro (K_1^{\ast} \times \ldots \times K_r^{\ast})
\underline{t}^{\underline{v}}$ has the natural structure of
semigroup.
\medskip

The extended semigroup $\widehat{S} \subset \mathds{Z}^r \times
(K_1^{\ast} \times \ldots \times K_r^{\ast})$ is invariant with
respect to the $K^{\ast}$-action. The factor space
\[
\pro \widehat{S} = \widehat{S} / K^{\ast}
\]
is called the projectivisation of the extended
semigroup (it is
also a graded semigroup in a natural sense).
\medskip

By the previous definitions, the projectivisation of the
extended semigroup can be described as a sum
\begin{equation} \label{eqn:pros}
\pro \widehat{S} = \sum_{\underline{v} \in \mathds{Z}^r} \pro
F_{\underline{v}} \cdot \underline{t}^{\underline{v}},
\end{equation}
where $\pro F_{\underline{v}}= F_{\underline{v}} / K^{\ast}$ is
the projectivisation of the fibre $F_{\underline{v}}$. For
$\underline{v} \in \widehat{S}$, the space $\pro
F_{\underline{v}}$ is the complement to an arrangement of
projective hyperplanes in a $\left (\dim_{K} \left (
J(\underline{v}) / J(\underline{v}+\underline{1}) \right ) - 1
\right)$-dimensional projective space $\pro \left (
J(\underline{v})/J(\underline{v}+\underline{1}) \right )$.

\section{The Poincar\'e series associated with the ring $\oo$}

\subsection{Multi-index filtrations and Poincar\'e series} \label{subsection:31}

Let us consider the multi-index filtration defined by the ideals
$\{J(\underline{v}) \}$,
for a given vector $\underline{v}=(v_1, \ldots , v_r) \in \mathds{Z}_{\ge 0}^r$. Since the ring
$\oo$ is a one--dimensional reduced ring, it is
Cohen-Macaulay, so the property of the ideals $J(\underline{v})$
of $\oo$ of containing regular elements is equivalent to the
property $\mm^{N} \supset J(\underline{v})$ for some integer
$N$. This means that the filtration $\{ J(\underline{v}) \}$ is in this case
finitely determined, and therefore every subspace $J(\underline{v})$  has finite codimension $\ell (\underline{v})$ in
 $\oo$. In particular, the dimension
 \[
 c(\underline{v}):= \dim_{k} \left (J(\underline{v})/J(\underline{v}+\underline{1})
\right )
 \]
 is also finite, because of the relation
 \[
 \ell(\underline{v}+\underline{1})=\ell(\underline{v})+c(\underline{v}).
 \]

This multi-index filtration can also be described in terms of a
multi-variable Laurent series

\begin{eqnarray} \label{eqn:lt}
L (t_1, \ldots , t_r) & := & \sum_{\underline{v} \in \mathds{Z}^r} c(\underline{v})
 \cdot \underline{t}^{\underline{v}},
\end{eqnarray}
where
$\underline{t}^{\underline{v}}:=t_1^{v_1} \cdot \ldots \cdot
t_r^{v_r}$. We will also write $L(\underline{t})$ instead of
$L(t_1, \ldots , t_r)$ if the number of variables is clear from
the context.  The series (\ref{eqn:lt})  is an element of the $\mathds{Z}[t_1, \ldots , t_r]$-module
(or even a module over the ring $\mathds{Z}[t_1, \ldots , t_r, t_1^{-1},
\ldots , t_r^{-1}]$) $\mathds{Z}[\![t_1, \ldots , t_r, t_1^{-1},
\ldots , t_r^{-1}]\!]$ (but not a ring, although the polynomial ring
$\mathds{Z}[t_1, \ldots , t_r]$ can be in a natural way considered
to be embedded into it). That is, $L (\underline{t})$ is a Laurent
series infinitely long in all directions (given by the vectors of the standard basis 
$\{\underline{1}_{1}, \ldots , \underline{1}_{r} \}$) since $c(\underline{v}):=\dim_{\mathds{C}}
\left ( J(\underline{v})/J(\underline{v}+\underline{1}) \right )$
can be positive for $\underline{v}$ with (some) negative
components $v_i$ as well. It is easy to see that, along each line
in the lattice $\mathds{Z}^r$ parallel to a coordinate one, the
coefficients $c(\underline{v})$ stabilize in each direction; that
is, if $v_i^{\prime}$ and $v_i^{\prime \prime}$ are negative, or
if $v_i^{\prime}$ and $v_i^{\prime \prime}$ are positive and large
enough, then $c(v_1, \ldots , v_i^{\prime}, \ldots v_r)=c(v_1,
\ldots , v_i^{\prime \prime}, \ldots v_r)$. This implies that
\begin{equation} \label{eqn:Pseries}
P^{\prime}(t_1, \ldots t_r ) = \prod_{i=1}^{r} (t_i -1) L(t_1, \ldots , t_r)
\end{equation}
is a power series in $\underline{t} = t_1\cdot \ldots \cdot t_r$ (cf. \cite{CDG-Alex}). 
\medskip

If we write $P^{\prime}(t_1, \ldots, t_r)= \sum_{\underline{v} \in \ZZ^r} p^{\prime} (\underline{v}) \underline{t}^{\underline{v}}$, then for every $\underline{v} \in \ZZ^r$ we have
\begin{eqnarray}
p^{\prime}(\underline{v}) & = & (-1)^r \sum_{j=0}^{r} (-1)^{j} \sum_{1 \le i_1 < \ldots < i_j \le r} c(\underline{v}-\underline{1}_{\{i_1, \ldots , i_j \}})  \nonumber \\
 & = & (-1)^ r \sum_{J \subseteq I} (-1)^{\sharp (J)} c(\underline{v}-\underline{1}_J).  \label{eqn:asterisco}
\end{eqnarray}

Notice that $P^{\prime} (t_1,\ldots , t_r)$ is a polynomial in the indeterminates $t_1, \ldots, t_r$ with coefficients in $\ZZ$. Moreover, if $p^{\prime}(\underline{v}) \ne 0$, then $\underline{0} \le \underline{v} \le \delta$, where $``\le"$ means the ordering given by the product.

\begin{Rem}
The dimension $c(\underline{v})$ depends only on the value semigroup $S$. In fact, if $K \subseteq I$, let
\[
\varnothing \subsetneq K_1 \subsetneq K_2 \subsetneq \ldots \subsetneq K_{h-1} \subsetneq K_h=I
\]
be such that $\Delta^{K_i}(\underline{v}) \ne \varnothing$ and $h$ is maximal with these properties. If we choose $g_i \in J(\underline{v})$ such that $\underline{v}(g_i) \in \Delta^{K_i}(\underline{v})$ for $0 \le i \le h$, then $\{j_{\underline{v}}(g_i) \mid 0 \le i \le h-1\}$ is a basis of $C(\underline{v})$, therefore $h=c(\underline{v})$ (cf. \cite{CDGext}).
\end{Rem}

\subsection{Computation of the coefficients of the Poincar\'e series} \label{Subsection:23}

Let $\underline{v} \in \ZZ^r$ and $i \in I$. The vector space $C(\underline{v},i)=J(\underline{v})/J(\underline{v}+\underline{1}_{\{i\}})$ can be identified with a vector subspace of the complex line. We will denote $c(\underline{v},i):=\dim_{k}C(\underline{v},i)$. It is clear that $c(\underline{v},i)=1$ if and only if there exists $\underline{w} \in S$ such that $w_i=v_i$ and $w_j \ge v_j$ for every $j \in I$. For any reordering $\{i_1, \ldots , i_r\}=\{1, \ldots , r\}$ of the set $I$ one has
\[
C(\underline{v}) \cong \bigoplus_{j=0}^{r-1} C(\underline{v}+\underline{1}_{\{i_1, \ldots , i_j \}}, i_{j+1}),
\]
and therefore
\[
c(\underline{v})= \sum_{j=0}^{r-1}c(\underline{v}+\underline{1}_{\{i_1, \ldots , i_j \}}, i_{j+1}).
\]

Let $\underline{v}, \underline{w} \in \ZZ^r$ and $\underline{v} \le \underline{w}$; for every $i \in I$ there exists a canonical linear map $\varphi_{\underline{w},\underline{v},i}: C(\underline{w},i) \rightarrow C(\underline{v},i)$ induced by the inclusion $J(\underline{w}) \to J(\underline{v})$. If, moreover, $w_i=v_i$, then the map $\varphi_{\underline{w},\underline{v},i}$ is a monomorphism, i.e., $0 \le c(\underline{w},i) \le  c(\underline{v},i) \le 1$.
\medskip

Let $\underline{v} \in \ZZ^r$ and fix $i \in I$. Define
\begin{eqnarray}
p_i(\underline{v}) & := & (-1)^{r-1} \sum_{J \subseteq I \setminus \{i\}} (-1)^{\sharp (J)} c(\underline{v}+\underline{1}-\underline{1}_{\{i\}}-\underline{1}_J, i) \nonumber \\
 & = & (-1)^r \sum_{i \in J^{\prime} \subseteq I} (-1)^{\sharp (J^{\prime})} c(\underline{v}+\underline{1}-\underline{1}_{J^{\prime}},i) \nonumber
\end{eqnarray}
and let
\[
P_i(t_1, \ldots , t_r) = \sum_{\underline{v} \in \ZZ^r} p_i (\underline{v}) \underline{t}^{\underline{v}}.
\]

The following result may be proved in much the same way as  Proposition 8 in \cite[p.~1649]{DGN}, we include however the proof by the sack of completeness.

\begin{Theo}
Let $\underline{v} \in \ZZ^r$ and $i \in I$. Then $p^{\prime}(\underline{v})=-p_i(\underline{v})+p_i(\underline{v}-\underline{1})$ and therefore
\[
(t_1 \cdot \ldots \cdot t_r -1) \cdot P_i(t_1, \ldots , t_r)= P^{\prime}(t_1, \ldots , t_r).
\]
As a consequence $P_i(t_1, \ldots, t_r)$ does not depend on $i$ (it will be denoted $P(t_1, \ldots , t_r)$ in the sequel).
\end{Theo}

\dem~The coefficient $p^{\prime}(\underline{v})$ can be written as
\[
p^{\prime}(\underline{v})= (-1)^r \sum_{i \notin J \subseteq I} (-1)^{\sharp (J)} (c(\underline{v}-\underline{1}_J)-c(\underline{v}-\underline{1}_{\{i\}}-\underline{1}_J)).
\]
Let $\{i_1, \ldots , i_{r-1}\}$ be any subset of the set $I \setminus \{i\}$. Then one can compute each summand above as
\[
c(\underline{v}-\underline{1}_J)= \sum_{j=0}^{r-2} c(\underline{v}-\underline{1}_J+\underline{1}_{\{i_1, \ldots , i_j \}}, i_{j+1})+c(\underline{v}-\underline{1}_J+\underline{1}_{I \setminus \{i\}},i)
\]
\[
c(\underline{v}-\underline{1}_J-\underline{1}_{\{i\}}) = c(\underline{v}-\underline{1}_J-\underline{1}_{\{i\}},i) + \sum_{j=0}^{r-2} c(\underline{v}-\underline{1}_J+\underline{1}_{\{i_1, \ldots , i_j \}},i_{j+1})
\]
and as a consequence one has
\begin{eqnarray}
p^{\prime}(\underline{v}) & = & (-1)^r \sum_{i \notin J \subseteq I} (-1)^{\sharp (J)} (c(\underline{v}-\underline{1}_J-\underline{1}_{\{i\}}+\underline{1},i)-c(\underline{v}-\underline{1}_{\{i\}}-\underline{1}_J,i)) \nonumber \\
& = & -p_i(\underline{v}) + p_i (\underline{v}-\underline{1}). \nonumber
\end{eqnarray}
\qed

\begin{Defi}
The  polynomial $P(t_1, \ldots, t_r)$ is called the Poincar\'e polynomial associated with $\oo$.
\end{Defi}

The dimensions $c(\underline{v})$ satisfy the following symmetry property (cf. \cite[Corollary 3.7]{CDK}; remember that we assumed the ring $\oo$ to be residually rational):

\begin{Theo}[Campillo, Delgado, Kiyek] \label{theo:cdk1}
Let $\underline{v} \in \mathbb{Z}^r$. Then
\[
c(\underline{v})+c(\delta- \underline{v}- \underline{1}) \le r.
\]
Furthermore, the ring $\oo$ is Gorenstein if and only if the equality holds.
\end{Theo}

The Poincar\'e polynomial associated with $\oo$ satisfies the following functional equation if the ring $\oo$ is assumed to be Gorenstein:

\begin{Theo}[Campillo, Delgado, Kiyek] \label{theo:cdk2}
If $\mathcal{O}$ is Gorenstein, then we have
\[
P^{\prime} (t_1, \ldots , t_r) + (-1)^{r} \underline{t}^{\delta} P^{\prime}(t_1^{-1}, \ldots , t_r^{-1})=0.
\]
\end{Theo} 
This formula gives immediately the relation 
\[
P(t_1, \ldots, t_r)=(-1)^{r} \underline{t}^{\tau} P(t_1^{-1}, \ldots , t_r^{-1}),
\]
with $\tau:=\delta-\underline{1}$, hence the coefficients satisfy
\[
p(\underline{v}) = (-1)^r p(\tau - \underline{v})
\] 
for all $\underline{v} \in \ZZ^r_{\ge 0}$. Next lemma gives some computations on the coefficients $p(\underline{v})$ of the Poincar\'e series $P(t_1, \ldots , t_r)$ of $\oo$.

\begin{Prop} \label{Lem:1}
We have:
\begin{itemize}

\item[(1)] If $\underline{v} \notin S$, then $p(\underline{v})=0$.

\item[(2)] If $\underline{v} \in S$ is not a maximal element of $S$, then $p(\underline{v})=0$.

\item[(3)] If $\underline{v}$ is an absolute maximal, then $p(\underline{v})=1$; if $\underline{v}$ is a relative maximal, then $p(\underline{v})=(-1)^{r}$.

\end{itemize}
\end{Prop}

\proof

We can fix any index $i \in I$ in order to compute $p(\underline{v})=p_i(\underline{v})$; then let us fix without loss of generality $i=r$. If $\{i_1, \ldots, i_{r-1}\}$ is any subset of $\{1, \ldots, r-1\}$, we have
\begin{eqnarray}
0  & \le &  c(\underline{v}+\underline{1}-\underline{1}_{\{r\}},r) \le c(\underline{v}+\underline{1}-\underline{1}_{\{r\}}-\underline{1}_{\{i_1\}},r) \le \ldots \nonumber \\
  & \ldots \le & c(\underline{v}+\underline{1}-\underline{1}_{\{r\}} -\underline{1}_{\{i_1, \ldots , i_{r-1}\}},r)=c(\underline{v},r) \le 1. \nonumber
\end{eqnarray}
If $c(\underline{v}+\underline{1}-\underline{1}_{\{r\}},r)=1$, then all terms involved in the expression
\[
p_r(\underline{v})=(-1)^{r-1} \sum_{j=0}^{r-1} (-1)^j \sum_{1 \le i_1 < \ldots < i_j \le r-1} c(\underline{v}+\underline{1}-\underline{1}_{\{i_1, \ldots , i_{j}\}},r)
\]
are equal to $1$. As a consequence we have
\[
p_r(\underline{v})=(-1)^{r-1} \sum_{j=0}^{r-1} (-1)^j {r-1 \choose j} = (-1)^{r-1}(1-1)^{r-1}=0. 
\]
If $c(\underline{v},r)=0$, then all terms are equal to $0$ and therefore $p_r(\underline{v})=0$.
Now, let $\underline{v} \notin S$; then there exists $i \in \{1, \ldots ,r\}$ such that $c(\underline{v},i)=0$ and by the previous computations we have $p(\underline{v})=p_i(\underline{v})=0$. If $\underline{v} \in S$ is not maximal, then there exists an index $i$ such that $\Delta_{\{i\}}(\underline{v}) \ne \varnothing$, but this is equivalent to say $c(\underline{v}+\underline{1}-\underline{1}_{\{i\}},i)=1$, and therefore $p(\underline{v})=p_i(\underline{v})=0$. Let $\underline{v} \in S$ be an absolute maximal; then $\Delta_K(\underline{v}) \ne \varnothing$ for every $K \subset I$, and $c(\underline{v},r)=1$ and $c(\underline{v}+\underline{1}-\underline{1}_{\{r\}}-\underline{1}_J,r)=0$ for every $J \subset I \setminus \{r\}$. Hence $p(\underline{v})=(-1)^{r-1}(-1)^{r-1}=1$. The condition for $\underline{v}$ to be a relative maximal is that $\Delta_{K}(\underline{v}) \ne \varnothing$ for every $K$ such that $\sharp (K) \ge 2$ (and it is empty if $\sharp (K)=1$). Therefore $c(\underline{v}+\underline{1}-\underline{1}_{\{r\}},r)=0$ and $c(\underline{v}+\underline{1}-\underline{1}_{\{r\}}-\underline{1}_J,r)=1$ for every $\varnothing \ne J \subset I \setminus \{r\}$. Then we have
\[
p(\underline{v})=(-1)^{r-1} \sum_{j=1}^{r-1} (-1)^j {r-1 \choose j} = (-1)^r.
\]
\qed

\begin{Rem}
From the above computations we get nice descriptions of $P(t_1, \ldots , t_r)$ for $r \in \{2,3\}$. Indeed, if $r=2$, there is no difference between relative and absolute maximals and we have
\[
P(t_1,t_2) = \sum_{\underline{v} \mathrm{~maximal}} \underline{t}^{\underline{v}}.
\]
If $r=3$, then we have only relative and absolute maximals and we obtain
\[
P(t_1,t_2,t_3)= -\sum_{\underline{v} \mathrm{~rel.~max.}} \underline{t}^{\underline{v}} +  \sum_{\underline{v} \mathrm{~abs.~max.}} \underline{t}^{\underline{v}}
\]
\end{Rem}

\section{Two possible definitions and its equivalence}

The definition of Poincar\'e series for $\mathcal{O}$ that we have chosen as a ``tentative" is not the only possible. First of all, if a vector $\underline{v} \in \mathds{Z}^r$ belongs to $S$, then it seems to be natural to take at least the  ``dimension" of its fiber $F_{\underline{v}}$ in the extended semigroup as a way to measure (or ``count'') the elements $g \in \mathcal{O}$ with $\underline{v}(g)=\underline{v}$. However, $F_{\underline{v}}$ is not a vector subspace of $k^r$ but a linear subspace minus some hyperplanes, more precisely $F_{\underline{v}}=C(\underline{v}) \cap (k \setminus \{0\})$, also $\overline{F_{\underline{v}}}=C(\underline{v})$. As a consequence, any definition of Poincar\'e series must take as a coefficient of $\underline{t}^{\underline{v}}$ at least the dimension of $C(\underline{v})$. Note that for any $\underline{v} \in \mathds{N}^r$ (indeed for any $\underline{v} \in \mathds{Z}^r$ with at least one coordinate greater that or equal to zero) the vector space $C(\underline{v})$ is different from $0$, even in the case $F_{\underline{v}}=\varnothing$ if $\underline{v} \notin S$.
\medskip

Another natural definition for the Poincar\'e series is
\[
\widetilde{L}_S(\underline{t})=\sum_{\underline{v} \in \mathds{N}^r} c(\underline{v}) \underline{t}^{\underline{v}} \in \mathds{Z}[\![t_1, \ldots , t_r]\!].
\] 
Using similar computations as those in the subsection \ref{Subsection:23} we can also take the polynomial
\[
\widetilde{P}_S(\underline{t})=(t_1-1) \cdot \ldots \cdot (t_r -1) \widetilde{L}_S(\underline{t}) =\sum_{\underline{v}\in \mathds{N}^r} \widetilde{p}(\underline{v})\underline{t}^{\underline{v}}.
\]
For any $\underline{v}$ with $v_i>0$ for $i \in I$ it obviously holds $\widetilde{p}(\underline{v})=p^{\prime}(\underline{v})$, hence the differences between $\widetilde{P}$ and $P^{\prime}$ appear only in the intersections of $\mathds{N}^r$ with the coordinate hyperplanes. More precisely, let $A \subseteq I$, and let $S_A$ be the semigroup arising from $S$ by projection on the set of indices $A$ (i.e. $S_A$ is the value semigroup associated with the curve consisting of the branches given by $A$). We have:

\begin{Theo} \label{Theo:uno}
If $r=1$, then $\widetilde{P}_S = P^{\prime}_S$. For $r \ge 2$ we have
\[
\widetilde{P}_S(\underline{t})=\sum_{\varnothing \ne A \subseteq I} (-1)^{\sharp (I \setminus A)} P^{\prime}_{S_A}(\underline{t}_A),
\]
where $\underline{t}_A:=\prod_{i \in A} t_i$. As a consequence, one has
\[
P^{\prime}_S(\underline{t})=\sum_{\varnothing \ne A \subseteq I} (-1)^{\sharp (I \setminus A)} \widetilde{P}_{S_A}(\underline{t}_A).
\]
\end{Theo}

\dem~Take $J \subseteq I$, and consider $\Lambda_J (\underline{t}):=\sum_{\underline{v}\in \overline{\Delta}^J(\underline{0})} \widetilde{p}(\underline{v}) \underline{t}^{\underline{v}}$ so that
\[
\widetilde{P}_S (\underline{t})=\sum_{J \subseteq I} \Lambda_J (\underline{t}_J).
\]
Let $J \ne \varnothing$. It is easily seen that
\[
\Lambda_J (\underline{t}_J)=(-1)^{\sharp (I \setminus J)} (P^{\prime}_{S_J}(\underline{t}_J)-p^{\prime}_{s_J}(\underline{v})).
\]
We know that $\widetilde{p}_S(\underline{v})=p(\underline{v})=(-1)^r$ and that $c(\underline{v}-\underline{1}_J)=1$ for $J \subseteq I$, $J \ne I$, and $c(\underline{v}-\underline{1}_J)=0$ for $J = I$, therefore
\[
p^{\prime}_S(\underline{v})=p^{\prime}(\underline{v})=(-1)^r \sum_{J \subsetneq I} =(-1)^r \sum_{i=0}^{r-1} (-1)^{i} {r \choose i} = (-1)^{2r-1}=-1,
\]
and so
\[
\Lambda_J(\underline{t}_J) =(-1)^{\sharp(I \setminus J)} (P^{\prime}_{S_J}(\underline{t}_J)+1).
\]
We have
\[
\widetilde{P}_S(\underline{t}) = \sum_{\varnothing \ne J \subseteq I} \Lambda_J (\underline{t}_J)+ \Lambda_{\varnothing} (\underline{t}_{\varnothing}),
\]
where $\Lambda_{\varnothing} (\underline{t}_{\varnothing})=\widetilde{p}(\underline{0})=(-1)^r$, hence we obtain:
\begin{eqnarray}
\widetilde{P}_S(\underline{t}) & = & \sum_{\varnothing \ne J \subseteq I} (-1)^{\sharp(I \setminus J)} (P^{\prime}_{S_J}(\underline{t}_J)-p^{\prime}_{s_J}(\underline{v}))+(-1)^r \nonumber \\
& = & \sum_{\varnothing \ne J \subseteq I} (-1)^{\sharp(I \setminus J)} P^{\prime}_{S_J}(\underline{t}_J)  -  \sum_{\varnothing \ne J \subseteq I} (-1)^{\sharp(I \setminus J)} p^{\prime}_{s_J}(\underline{v})   +(-1)^r \nonumber \\
& = & \sum_{\varnothing \ne J \subseteq I} (-1)^{\sharp(I \setminus J)} P^{\prime}_{S_J}(\underline{t}_J)  -  \sum_{\varnothing \ne J \subseteq I} (-1)^{\sharp(I \setminus J)} (-1)  +(-1)^r \nonumber \\
& = & \sum_{\varnothing \ne J \subseteq I} (-1)^{\sharp(I \setminus J)} P^{\prime}_{S_J}(\underline{t}_J)  +  \sum_{\varnothing \ne J \subseteq I} (-1)^{\sharp(I \setminus J)}  +(-1)^r \nonumber \\
& = & \sum_{\varnothing \ne J \subseteq I} (-1)^{\sharp(I \setminus J)} P^{\prime}_{S_J}(\underline{t}_J)  +  \sum_{i=1}^r (-1)^{i} {r \choose i} +(-1)^r \nonumber \\
& = & \sum_{\varnothing \ne J \subseteq I} (-1)^{\sharp(I \setminus J)} P^{\prime}_{S_J}(\underline{t}_J)  - (-1)^r +(-1)^r \nonumber \\
& = & \sum_{\varnothing \ne J \subseteq I} (-1)^{\sharp(I \setminus J)} P^{\prime}_{S_J}(\underline{t}_J), \nonumber
\end{eqnarray}
as desired. \qed

\begin{Rem}
Formulae above express the equivalence between the sets of polynomials $\{ P^{\prime}_{S_A} \mid A \subseteq I\}$ and $\{\widetilde{P}_{S_A} \mid A \subseteq I\}$. In fact, $\widetilde{P}_{S_A}$ can be obtained from $\widetilde{P}_S$ just by substituting $t_i=0$ for any $i \notin A$ and multiplying by $(-1)^{r-\sharp(A)}$. Then the polynomial $\widetilde{P}_S$ is equivalent to the set of polynomials $\{ P^{\prime}_{S_A} \mid A \subseteq I\}$. From this point of view there is no difference between the two definitions (see also Section \ref{section:5}).
\end{Rem}

For some purposes it is sometimes useful to consider codimensions of the ideals $J(\underline{v})$ as coefficients of the corresponding Poincar\'e series, what leads to a third definition:
\[
H(\underline{t}):=\sum_{\underline{v} \in \mathbb{N}^r}\ell (\underline{v}) \underline{t}^{v},
\]
where $\ell (\underline{v})= \dim (\oo/J(\underline{v}))$ is the codimension of the ideal $J(\underline{v})$ in the ring $\oo$ (cf. Subsection \ref{subsection:31}). As a consequence of Theorem \ref{Theo:uno} one can prove the following relation between the Poincar\'e series $H(\underline{t})$ and $P(\underline{t})$:

\begin{Cor} \label{cor:H}
We have
\[
\prod_{i=1}^{r} (1-t_i) H(\underline{t})= \sum_{\varnothing \ne A \subseteq I} (-1)^{\sharp (A)-1} \underline{t}_{A}P_{S_A}(\underline{t}_A).
\]
\end{Cor}

\dem~Set $Q(\underline{t}):=\underline{t} \cdot P(\underline{t})$. It is enough to realise that
\[
Q(\underline{t})=(-1)^r \prod_{i=1}^r(1-t_i) H(\underline{t})
\]
(using the fact that $\ell (\underline{v} + \underline{1})=\ell (\underline{v}) + c(\underline{v})$ for any $\underline{v} \in \mathds{Z}^r$).
By applying Theorem \ref{Theo:uno} to the series $Q(\underline{t})$, we are done. \qed

\section{Behaviour by projections} \label{section:5}

In this section we want to analyse the behaviour by projections of the multi-variable Poincar\'e series, i.e., geometrically speaking, the behaviour of the Poincar\'e series of a reduced but not irreducible plane curve singularity by removing some of its branches. Consider
\[
 J=\{i_1, \ldots, i_h \} \subset I=\{1, \ldots,r \}.
\] 
Denoting
by $\mathrm{pr}_J$ the projection over the indexes $J$ of
$\mathbb{Z}_{+}^{r}$ in $\mathbb{Z}_{+}^{\sharp J}$. We set
\[
S_{\{i_1 , \ldots, i_h\}}=S_J:=\mathrm{pr}_J (S)= \{\mathrm{pr}_J
(\underline{v}) \mid \underline{v} \in S \}.
\]
That is, the semigroup $S_J$ corresponds to the projection of
$S_I$ on the branches indexed by $J$.
\medskip

Let $\xi_{i,j}$ denote the intersection multiplicity between the branches $C_i$ and $C_j$ given by the local rings $V_i$ and $V_j$, for every $i,j \in \{1, \ldots ,r\}$ with $i \ne j$.

\begin{Theo} \label{Theo:proj}
Let $P_{S_{\{r\}}}(t_1, \ldots , t_{r-1})$ be the Poincar\'e polynomial associated with the ring consisting of the branches $C_1, \ldots , C_{r-1}$ corresponding to the local rings $V_1, \ldots , V_{r-1}$, and $P(t_1, \ldots , t_r)$ the Poincar\'e polynomial of $\oo$. We have
\[
P(t_1, \ldots , t_{r-1},1)=(1-t_1^{\xi_{1,r}} \cdot \ldots \cdot t_{r-1}^{\xi_{r-1,r}}) P_{S_{\{r\}}}(t_1, \ldots, t_{r-1}).
\]
The same formula holds for $P^{\prime}$.
\end{Theo}

\dem~By the description of the Poincar\'e series $P^{\prime}_{S_{\{r\}}}(\underline{t})$ given by (\ref{eqn:asterisco}) we have
\begin{eqnarray}
(1-\underline{t}^{\xi^{(r)}}) P^{\prime}_{S_{\{r\}}}(t_1, \ldots , t_{r-1}) & = & \sum_{\underline{v} \in S_{\{r\}}} p^{\prime}_{S_{\{r\}}} (\underline{v}) \underline{t}_{r}^{\underline{v}} - \sum_{\underline{v} \in S_{\{r\}}} p^{\prime}_{S_{\{r\}}}(\underline{v}) \underline{t}_{r}^{\underline{v}} \underline{t}_{r}^{\xi^{(r)}} \nonumber \\ 
&=&  \sum_{\underline{v} \in S_{\{r\}}} A(\underline{v}) \underline{t}_{r}^{\underline{v}} \nonumber
\end{eqnarray}
for $A(\underline{v}):=p^{\prime}_{S_{\{r\}}}(\underline{v})-p^{\prime}_{S_{\{r\}}} (\underline{v}-\xi^{(r)} )$ and with $\underline{t}_{r}^{\underline{v}}:=t_1^{v_1} \cdot \ldots \cdot t_{r-1}^{v_{r-1}}$ and $\underline{t}_{r}^{\xi^{(r)}}:=t_1^{\xi_{1,r}} \cdot \ldots \cdot t_{r-1}^{\xi_{r-1,r}}$. The key point here is to understand the behaviour of the coefficients of the Poincar\'e series under semigroup projections. Let us present this for $S_{\{r\}}$, and the result follows by induction.
\medskip

If we consider coefficients $p_{S_{\{r\}}}^{\prime}(\underline{w})$ in the fiber of such a semigroup projection of a fixed element $\underline{w}=(w_1, \ldots, w_{r-1}) \in S_{\{r\}}$, we realise first that the coefficients $p^{\prime}_{S_{\{r\}}} (\underline{w})$ of the projection series  can be written as a finite sum $\sum_{\nu=0}^{u} p^{\prime}(\underline{w},\nu)$. Indeed, our aim is to show the equality
\[
A(\underline{v})=\sum_{\nu=0}^{u} p^{\prime}(\underline{w},\nu),
\]
and then we will be done. A first remarkable fact is the following: From the definitions of $P^{\prime}_S(t_1, \ldots , t_r)$ and $P^{\prime}_{S_{\{r\}}}(t_1, \ldots , t_{r-1})$, the only remaining summands in the coefficients $p_S^{\prime}(\underline{v})$ are
\[
(-1)^r \sum_{r \notin J} (-1)^{\sharp (J)} c((\underline{w},u)-\underline{1}_{J}) + (-1)^r \sum_{r \in J} (-1)^{\sharp (J)} c((\underline{w},1)-\underline{1}_{J}).
\]
Next we apply the symmetry properties given by Theorem \ref{theo:cdk1} to show
\[
c((\underline{w},u)-\underline{1}_{J})=r-c(\delta-\underline{1}-(\underline{w},u)+\underline{1}_{J}),
\]
because our curve is assumed to be plane, hence Gorenstein. Denote $J^{\prime}:=I \setminus J$. Since $c((\underline{w},u)-\underline{1}_{J})=r-c((\delta-\underline{w},0)-\underline{1}_{J^{\prime}})$, the term $(-1)^r \sum_{r \notin J} (-1)^{\sharp (J)} c((\underline{w},u)-\underline{1}_{J}) $ is equal to
\begin{align*}
& (-1)^r \sum_{r \notin J} (-1)^{\sharp (J)} (r-c((\delta-\underline{v},0)-\underline{1}_{I \setminus J}))   =\\
 =& (-1)^r \sum_{r \notin J} (-1)^{\sharp (J)} - (-1)^r \sum_{r \notin J} (-1)^{\sharp (J)} c((\delta-\underline{w},0)-\underline{1}_{J^{\prime}})  = (\ast) 
\end{align*}
One may easily prove that $\sum_{r \notin J} (-1)^{\sharp (J)}=\sum_{i=0}^{r-1} (-1)^i {r-1 \choose i}=0$, hence 
\begin{eqnarray}
(\ast) & = & (-1)^{r-1} \sum_{r \in J^{\prime}} (-1)^{r-\sharp (J^{\prime})} c ((\delta-\underline{w},0)-\underline{1}_{J^{\prime}}) \nonumber \\
&=& (-1)^r \Big ((-1)^{r-1} \sum_{r \in J^{\prime}} (-1)^{\sharp (J^{\prime})} c ((\delta-\underline{w},0)-\underline{1}_{J^{\prime}}) \Big ). \nonumber
\end{eqnarray}
On the other hand, it is easily checked that
\[
(-1)^r \sum_{r \in J^{\prime}} (-1)^{\sharp (J^{\prime})} c ((\underline{w},h)-\underline{1}_{J^{\prime}})=p^{\prime}_{S_{\{r\}}}(\underline{w})
\]
for $h \in \{0,1\}$, and therefore
\[
\sum_{\nu=1}^{u} p^{\prime} (\underline{w},\nu)=(-1)^{r-1}p^{\prime}_{S_{\{r\}}}(\delta-\underline{w}) + p^{\prime}_{S_{\{r\}}}(\underline{w}).
\]
By applying Theorem \ref{theo:cdk2} to the coefficients $p^{\prime}_{S_{\{r\}}}(\underline{w})$ we get
\[
(-1)^{r-1}p^{\prime}_{S_{\{r\}}}(\delta-\underline{w})=-p^{\prime}_{S_{\{r\}}}(\underline{w}-\xi^{(r)})
\]
(since $\delta-\underline{w}=\delta(S_{\{r\}})+(\xi_{1,r}, \ldots ,\xi_{r-1,r})-\underline{w}=\delta (S_{\{r\}})-(\underline{w}-\xi^{(r)})$). Therefore
\[
\sum_{\nu=1}^{u} p^{\prime} (\underline{w},\nu)=p^{\prime}_{S_{\{r\}}}(\underline{w}) -p^{\prime}_{S_{\{r\}}}(\underline{w}-\xi^{(r)})= A(\underline{w}),
\]
as desired.
\qed

\begin{Cor}
More generally, if $A \subseteq I$, then we have
\[
P_S(\underline{t})\mid_{\{t_i=1 \mid i \notin A\}} = \prod_{j \notin A} (1-\underline{t}_A^{\xi_{j_A}}) P_{S_A}(\underline{t}_A),
\]
where $\underline{t}_A^{\xi_{j_A}}:=\prod_{i \in A} t_i^{\xi_{i,j}}$ and $P_{S_A}(\underline{t}_A)$ is the Poincar\'e polynomial of $S_A$.
\end{Cor}

Next result seems to involve non-trivial computations on $P$:

\begin{Theo}
The intersection multiplicities $\xi_{i,j}$ between couples of branches of $\oo$ can be computed from $P_S$.
\end{Theo}

\dem~First, we may apply recursively the projection formula to obtain
\[
P^{\prime}_S (1, \ldots 1,\underset{\underset{i}{\uparrow}}{t}, 1, \ldots, 1)=\prod_{\stackrel{j=1}{j \ne i} }^{r} (1-\underline{t}^{\xi_{i,j}}) \cdot P^{\prime}_{S_{\{i\}}}(t).
\]

At this point, we need a theorem of Campillo, Delgado and Gusein-Zade (cf.~\cite{CDG-i}; also~\cite{MK}) which expresses the Poincar\'e series of an irreducible plane curve singularity in terms of a minimal set of generators of the corresponding value semigroup. More precisely, assuming the set $\{\rho^{i}_0, \rho^{i}_1, \ldots , \rho^{i}_g\}$ to be a minimal set of generators of $S_{\{i\}}$, and writing $\theta^{i}_0:=\rho^{i}_0$, $\theta^{i}_{k+1}:= \mathrm{gcd}(\rho^{i}_0, \ldots , \rho^{i}_k)$ for $k \in \{1, \ldots , g\}$, and  $N^{i}_k = \frac{\theta^{i}_k}{\theta^{i}_{k+1}}$ for all $k \in \{1, \ldots , g\}$ and all $i\in \{1, \ldots , r\}$, we have
\[
P^{\prime}_{S_{\{i\}}}(t)=(t-1) \sum_{n \in S_{\{i\}}} t^n = (t-1) \cdot \frac{1}{1-t^{\rho^{i}_0}} \cdot \prod_{k=1}^{g} \frac{(1-t^{N^{i}_k \rho^{i}_k})}{(1-t^{\rho^{i}_k})}
\]
is a polynomial non divisible by $t-1$ (see e.~g.~\cite{CDG-Alex}), and because of the fact that $(t-1) \sum_{i=0}^{m-1}t^{i}=(t^m-1)$, it is also not divisible by $t^m-1$ for any natural number $m$. Therefore we can write
\[
P^{\prime}_S (1, \ldots 1,\underset{\underset{i}{\uparrow}}{t}, 1, \ldots, 1)=\prod_{\stackrel{j=1}{j \ne i} }^{r} (1-\underline{t}^{z_{i,j}}) \cdot P^{\prime}_{S_{\{i\}}}(t).\\
\] 
with $z_{i,j} \in \mathbb{N}$, and this factorisation is uniquely determined. By comparing now the factors $\prod_{\stackrel{j=1}{j \ne i} }^{r} (1-\underline{t}^{z_{i,j}})$ and $\prod_{\stackrel{j=1}{j \ne i} }^{r} (1-\underline{t}^{\xi_{i,j}})$ for every $i \in \{1, \ldots, r\}$, we obtain the following system of linear equations in the indeterminates $z_{i,j}$:
\[
\left.
\begin{array}{rcl}
     z_{1,2}+z_{1,3} + z_{1,4}+ \ldots + z_{1,r} & = & \xi_{1,2}+\xi_{1,3} + \xi_{1,4}+ \ldots + \xi_{1,r} 
  \\ z_{2,1}+z_{2,3} + z_{2,4}+ \ldots + z_{2,r}  & = & \xi_{2,1}+\xi_{2,3} + \xi_{2,4}+ \ldots + \xi_{2,r} 
  \\ z_{3,1}+z_{3,2} + z_{3,4}+ \ldots + z_{3,r}  & = & \xi_{3,1}+\xi_{3,2} + \xi_{3,4}+ \ldots + \xi_{3,r} 
  \\ & \vdots &
  \\z_{r,1}+z_{r,2} + z_{r,3}+ \ldots + z_{r,r-1}  & = & \xi_{r,1}+\xi_{r,2} + \xi_{r,3}+ \ldots + \xi_{r,r-1} 
\end{array}
\right\}
\]
The resolution of the system yields the desired result.
\qed

\begin{Rem}
As a corollary from the previous results one has that the polynomials $P_{S_A}$, for $A \subseteq I$, can be computed from $P_S$. This provides---for complex plane curve singularities---an alternative way to prove that the Alexander polynomial is equivalent to the topological type of the singularity, and it was already proved with topological methods in \cite{CDG-top}. In particular, they describe how to reconstruct the resolution tree from the Poincar\'e series. It is worth pointing out that Yamamoto was the first who proved that the Alexander polynomial determines the topology of an algebraic link by topological methods (\cite{Y}). Also, it provides a proof that $P^{\prime}_S$ and $\widetilde{P}_S$ are equivalent data.
\end{Rem}

\section{Generalised Poincar\'e series}

%The aim of this section is to show the relation between the generalised (or motivic) Poincar\'e series $P_g (t)$ and $\widehat{P}_g(t)$ via the results about projections already seen in the second section.

One can now pose the question Section 5 deals with in the context of Poincar\'e series of motivic nature (called generalised Poincar\'e series following the terminology introduced in \cite{cadegu11}). Let us present first the basic notions.
\medskip

Let $K_0(\Nu_{k})$ be the Grothendieck ring of quasi-projective
varieties.
It is generated by classes $[X]$ of such varieties subject to the
relations:
\begin{itemize}
    \item[i)] if $X_1 \cong X_2$ then $[X_1]=[X_2]$;
    \item[ii)] if $Y$ is Zariski-closed in $X$, then $[X]=[Y]+[X \setminus Y]$
\end{itemize}
(the multiplication is defined by the Cartesian product). Let
$\mathbb{L}$ be the class $[\mathbb{A}^1_{k}]$ of the complex
affine line. The class $\mathbb{L}$ is not equal to zero in the
ring $K_0 (\Nu_{k})$. Moreover the natural ring homomorphism
$\mathbb{Z}[X] \to K_0 (\Nu_{k})$ which sends $X$ to
$\mathbb{L}$ is an inclusion. Let $K_0 (\Nu_{k})_{(\mathbb{L})}$
be the localisation of the Grothendieck ring $K_0 (\Nu_{k})$ by
the class $\mathbb{L}$. The natural homomorphism
$\mathbb{Z}[X]_{(X)} \to K_0 (\Nu_{k})_{(\mathbb{L})}$ is an
inclusion as well.
\medskip

The generalised Euler characteristic $\chi_g (X)$ of a cylindric subset $X$
is the element $[Y] \cdot \mathbb{L}^{-d(p)}$ in the ring $K_0
(\Nu_{k})_{(\mathbb{L})}$, where $Y=\pi^{-1}_p(X)$ is a
constructible subset of $\mathbb{P}\oo$, as in the previous
subsection. Note that $\chi_g (X)$ is well-defined, because if
$X=\pi^{-1}_q(Y^{\prime})$, $Y^{\prime} \subset \pro J_{V,0}^q$
and $p \ge q$, then $Y$ is a locally trivial fibration over
$Y^{\prime}$ with the fibre $k^{d(p)-d(q)}$ and therefore $[Y] =
[Y^{\prime}] \cdot \mathbb{L}^{d(p)-d(q)}$.
\medskip

Let us take the usual Euler characteristic $\chi (\cdot)$.
Let $A$ (resp. $A^{\prime}$) be a subspace of $\oo$ of
finite codimension $a$ (resp. $a^{\prime}$) with $a^{\prime} >a$.
Let be $Q:=\mathbb{L}^{-1}$. Then one has
\begin{eqnarray}
\dim_{k} \left ( A/A^{\prime} \right) & = & \chi \left ( \pro A
\setminus \pro A^{\prime} \right ) \label{eqn:1} \\
 & = & \chi \left ( \pro \left (
A/ A^{\prime} \right ) \right ) . \label{eqn:2}
\end{eqnarray}
These equalities do not hold for the generalised Euler
characteristic. Ne\-vertheless we can take Equation (\ref{eqn:1})
as definition of a sort of ``generalised (or motivic) dimension"
just by taking (\ref{eqn:1}) with $\chi_g (\cdot)$ instead of the
usual $\chi (\cdot)$. Hence using the cellular decomposition of a
projective space one gets
\begin{eqnarray}
\chi_g (\pro A \setminus \pro A^{\prime}) & = & Q^a + Q^{a+1} +
\ldots + Q^{a^{\prime}-1} \nonumber \\
 & = & Q^{a+1} \cdot \frac{1-Q^{a^{\prime}-a}}{1-Q}. \label{eqn:Qlg}
\end{eqnarray}

If we set $\ell (\underline{v}):=\dim_{k} \left ( \oo /
J(\underline{v}) \right )$ and put $A=J(\underline{v})$,
$A^{\prime}=J(\underline{v}+1)$, $a=\ell (\underline{v})$ and
$a^{\prime}= \ell (\underline{v}+1)$ in Equation (\ref{eqn:Qlg}),
then we define a series
\[
L_g (\underline{t}):= \sum_{\underline{v} \in \mathds{Z}^r}
Q^{\ell (\underline{v})+1} \cdot \frac{1-Q^{\ell (\underline{v}+1)
- \ell (\underline{v})}}{1-Q} \cdot \underline{t}^{\underline{v}}
\]
which is a ``motivic version'' of the series $L(\underline{t})$
given by the formula (\ref{eqn:lt}). One can also see (cf.
\cite[Proposition 2]{cadegu11}):

\[
P_g (\underline{t}) = \frac{L_g(\underline{t}) \cdot
\prod_{i=1}^{r} (t_i-1)}{t_1 \cdot \ldots \cdot t_r-1}.
\]
On the other hand, the equality between Equations (\ref{eqn:1})
and (\ref{eqn:2}) does not hold when we replace $\chi (\cdot)$ with
$\chi_g (\cdot)$:
\[
\chi_g \left ( \pro A \setminus \pro A^{\prime}) \right ) \ne
\chi_g \left ( \pro \left ( A/ A^{\prime} \right ) \right ).
\]

For $A=J(\underline{v})$ and $A^{\prime}=J(\underline{v}+1)$
we get
\[
\chi_g \left ( \pro \left ( J(\underline{v})/J(\underline{v}+1)
\right ) \right ) =
\frac{Q^{\ell(\underline{v})+1}}{Q^{\ell(\underline{v}+
\underline{1})}} \cdot \frac{1-Q^{\ell
(\underline{v}+\underline{1})-\ell(\underline{v})}}{1-Q},
\]
therefore we define the corresponding series as
\[
\widehat{L}_g (\underline{t}):=\sum_{\underline{v} \in
\mathds{Z}^r}
\frac{Q^{\ell(\underline{v})+1}}{Q^{\ell(\underline{v}+
\underline{1})}} \cdot \frac{1-Q^{\ell
(\underline{v}+\underline{1})-\ell(\underline{v})}}{1-Q}
\underline{t}^{\underline{v}}.
\]

\begin{Lem}
\[
\widehat{L}_g (\underline{t})=\sum_{\underline{v} \in
\mathds{Z}^r}
\left [ \mathbb{P} (J(\underline{v})/J(\underline{v}+\underline{1})) \right ]
\underline{t}^{\underline{v}}.
\]
\end{Lem}

\dem~The class in $K_0(\Nu_{k})$ of the projectivisation of the vector space $J(\underline{v})/J(\underline{v}+\underline{1})$ (of finite dimension) is equal to $(\mathbb{L}-1)^{-1}(\mathbb{L}^{c(\underline{v})}-1)$. By setting $Q:=\mathbb{L}^{-1}$ this is equal to $(1-Q)^{-1}Q (Q^{\ell(\underline{v})-\ell(\underline{v}+\underline{1})}-1)$, and this coincides with the coefficient of the definition of $\widehat{L}_g (\underline{t})$. \qed
\medskip

The following result shows a ``motivic'' series analogous to the series resulting from taking Euler characteristic to the spaces $\pro
F_{\underline{v}}$ of the formula (\ref{eqn:pros}):
\begin{Prop} 
\[
\widehat{P}_g (\underline{t}) = \frac{ \widehat{L}_g
(\underline{t}) \cdot \prod_{i=1}^{r} (t_i-1)}{t_1 \cdot \ldots
\cdot t_r-1}.
\]
\end{Prop}

%One may also define an affine version of the series $\widehat{L}_g (\underline{t})$:
%\[
%\widehat{L}^{\mathrm{aff}}_g(\underline{t}):= \sum_{\underline{v} \in
%\mathds{Z}^r}
%\left [J(\underline{v})/J(\underline{v}+\underline{1})\right ]
%\underline{t}^{\underline{v}}.
%\]
%Notice that $\widehat{L}^{\mathrm{aff}}_g(\underline{t}):= \sum_{\underline{v} \in
%\mathds{Z}^r} \mathbb{L}^{c(\underline{v})} \underline{t}^{\underline{v}}$, and we set
%\[
%\widehat{P}^{\mathrm{aff}}_g (\underline{t}) := \frac{ \widehat{L}^{\mathrm{aff}}_g
%(\underline{t}) \cdot \prod_{i=1}^{r} (t_i-1)}{t_1 \cdot \ldots
%\cdot t_r-1}.
%\]
%Following the same arguments as in \cite{MZ} one has:
%\[
%\widehat{P}^{\mathrm{aff}}_g(\underline{t})=(\mathbb{L}-1) \widehat{P}_g(\underline{t}).
%\]
%\[
%\widehat{L}^{\mathrm{aff}}_g(\underline{t})=(\mathbb{L}-1) \widehat{L}_g(\underline{t}).
%\]
%\medskip
Notice that the series $\widehat{P}_g (\underline{t})$ coincides with $P(\underline{t})$ in the case $r=1$.
\medskip

As a difference from the classical Poincar\'e polynomial, the generalised Poincar\'e series $P_g(\underline{t})$ does not satisfy a property analogous to that of Theorem \ref{Theo:proj}, but it just forget components, as E. Gorsky has shown in \cite{Gor} for the complex case. Let us define
\[
\overline{P}_g(t_1, \ldots , t_r):=(1-Q t_1) \cdot \ldots \cdot (1-Qt_r) \cdot P_g(t_1, \ldots ,t_r).
\]

Moreover, consider the series $P_{g}^{S_{\{r\}}}(t_1, \ldots , t_{r-1})$ be the generalisd Poincar\'e series associated with the ring consisting of the branches $C_1, \ldots , C_{r-1}$ corresponding to the local rings $V_1, \ldots , V_{r-1}$, and 
\[
\overline{P}_{g}^{S_{\{r\}}}(t_1, \ldots , t_{r-1}):=(1-Q t_1) \cdot \ldots \cdot (1-Qt_{r-1}) \cdot  P_{g}^{S_{\{r\}}}(t_1, \ldots , t_{r-1})
\]
\begin{Prop} [Gorski]
For  $\oo$ reduced with $r>1$ branches we have
\[
\overline{P}_g(t_1, \ldots , t_{r-1},t_r=1)=(1-Q) \overline{P}_{g}^{S_{\{r\}}}(t_1, \ldots , t_{r-1}).
\]
If $\oo$ is irreducible, then $\overline{P}_g(t=1)=1$.
\end{Prop}

It implies in particular that these generalised Poincar\'e series are not affected by ``forgetting components'', as a difference with the non-motivic Poincar\'e polynomial (cf. Theorem \ref{Theo:proj}).
\medskip

As we did in the non-motivic case, one can also define a motivic version of the Poincar\'e series $H(\underline{t})$, namely
\[
H_g(\underline{t}):=\sum_{\underline{v} \in \mathds{Z}_{\ge 0}^r} [\mathcal{O}/J(\underline{v})] \underline{t}^{\underline{v}}.
\]
Similarly as done for the series $H(\underline{t})$ in Corollary \ref{cor:H} we may also prove the following result:
\begin{Prop} \label{prop:hg}
\[
\prod_{i=1}^{r} (1-t_i) H_g(\underline{t})= \sum_{\varnothing \ne A \subseteq I} (-1)^{\sharp (A)-1} \underline{t}_{A}P_{g_{S_A}}(\underline{t}_A).
\]
\end{Prop}

\section*{Acknowledgements}

The author wishes to express his gratitude to Prof. Dr. F\'elix Delgado de la Mata for stimulating conversations and helpful remarks about the topic.

\end{document}